\documentclass[onecolumn,authoryear]{els-mrw} 

\usepackage{amsmath,amssymb,amsfonts,amsthm,makeidx,graphicx}
\usepackage{txfonts}
\usepackage{helvet}

\usepackage{xspace}
\usepackage{subcaption}

\usepackage{accents}
\usepackage{xcolor}
\usepackage{pifont}

\newcommand{\D}{\mathbb{D}}
\newcommand{\N}{\mathbb{N}}
\newcommand{\R}{\mathbb{R}}

\newcommand{\e}{\boldsymbol{e}}
\newcommand{\x}{\boldsymbol{x}}
\newcommand{\y}{\boldsymbol{y}}
\newcommand{\w}{\boldsymbol{w}}
\newcommand{\uv}{\boldsymbol{u}}
\newcommand{\z}{\boldsymbol{z}}

\newcommand{\0}{\boldsymbol{0}}
\newcommand{\1}{\boldsymbol{1}}

\newcommand{\A}{\boldsymbol{A}}
\newcommand{\W}{\boldsymbol{W}}
\newcommand{\Ib}{\boldsymbol{I}}

\newcommand{\I}{\mathsf{I}}
\newcommand{\T}{\mathsf{T}}
\newcommand{\Lcal}{\mathcal{L}}
\newcommand{\Ecal}{\mathcal{E}}
\newcommand{\Ncal}{\mathcal{N}}

\DeclareMathOperator*{\argmin}{arg\,min}
\DeclareMathOperator{\prox}{prox}
\DeclareMathOperator{\refl}{refl}
\DeclareMathOperator{\proj}{proj}
\DeclareMathOperator{\fix}{fix}
\DeclareMathOperator{\blkdiag}{blk\,diag}

\newcommand{\norm}[1]{\left\lVert#1\right\rVert}

\newcommand{\km}{Krasnosel'ski\u{\i}-Mann\xspace}
\newcommand{\ccp}{\Gamma_0}

\newcommand{\ubar}[1]{\underaccent{\bar}{#1}}
\newcommand{\lmin}{\ubar{\lambda}}
\newcommand{\lmax}{\bar{\lambda}}
\newcommand{\comp}{\circ}

\newcommand{\cmark}{\textcolor{green!80!black}{\ding{51}}}
\newcommand{\xmark}{\textcolor{red}{\ding{55}}}

\begin{document}

\chapter{Multi-Agent Optimization and Learning: A Non-Expansive Operators Perspective}\label{chap1}

\author[1]{Nicola Bastianello}%
\author[2]{Luca Schenato}%
\author[2]{Ruggero Carli}%

\address[1]{\orgname{KTH Royal Institute of Technology}, \orgdiv{School of Electrical Engineering and Computer Science, and Digital Futures}, \orgaddress{Stockholm, Sweden}}
\address[2]{\orgname{University of Padova}, \orgdiv{Department of Information Engineering (DEI)}, \orgaddress{Padova, Italy}}

\maketitle

\begin{abstract}[Abstract]
Multi-agent systems are increasingly widespread in a range of application domains, with optimization and learning underpinning many of the tasks that arise in this context. Different approaches have been proposed to enable the cooperative solution of these optimization and learning problems, including first- and second-order methods, and dual (or Lagrangian) methods, all of which rely on consensus and message-passing. In this article we discuss these algorithms through the lens of non-expansive operator theory, providing a unifying perspective. We highlight the insights that this viewpoint delivers, and discuss how it can spark future original research.
\end{abstract}

\section{Introduction}\label{sec:introduction}
\begin{figure}[!ht]
\centering
\begin{subfigure}{.5\textwidth}
\centering
    \includegraphics[scale=0.75]{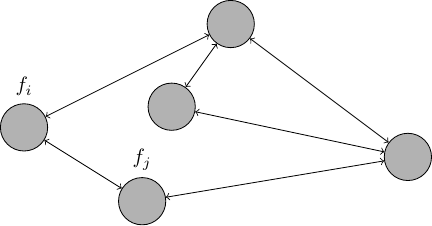}
    \caption{Distributed}
    \label{fig:architecture-distributed}
\end{subfigure}%
\begin{subfigure}{.5\textwidth}
\centering
    \includegraphics[scale=0.75]{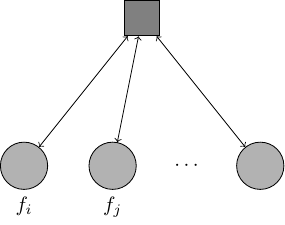}
    \caption{Federated}
    \label{fig:architecture-federated}
\end{subfigure}
\caption{Depictions of multi-agent optimization and learning architectures.}
\label{fig:architectures}
\end{figure}

The technological advancement of the past decade have led to a widespread deployment of multi-agent systems in a broad range of domains, including robotics, power grids, and traffic networks (\cite{nedic_distributed_2018}).
These systems are composed of agents equipped with computational and communication resources, interconnected with each other to pursue a cooperative task. There are two main communication architectures: \textit{distributed}, which employs peer-to-peer communications (cf. Figure~\ref{fig:architecture-distributed}), and \textit{federated}, in which the communications are mediated by a coordinator agent (cf. Figure~\ref{fig:architecture-federated}). In the following we focus on distributed scenarios, and briefly discuss federated ones in the conclusion.

A wide set of the tasks faced by multi-agent systems -- from coordination and collaborative decision-making to sensing and learning --  can be formulated as the distributed optimization problem:
\begin{equation}\label{eq:problem}
\begin{split}
    &\min_{x_i \in \R^n, \ i \in [N]} \sum_{i \in [N]} f_i(x_i) \\
    &\text{s.t.} \ x_1 = \ldots = x_N,
\end{split}
\end{equation}
where $[N] = \{ 1, \ldots, N \}$.
The objective in~\eqref{eq:problem} is for the agents to achieve consensus (as encoded by the constraints $x_1 = \ldots = x_N$) on a minimizer of the sum of local costs $\{ f_i \}_{i \in [N]}$ (which we denote as $x^* \in \argmin \sum_{i \in [N]} f_i(x)$).
The local cost functions are defined on data and measurements accessible to each of the agent; for example, in learning they are defined as empirical risk minimization costs
\begin{equation}\label{eq:erm}
    f_i(x) = \frac{1}{m_i} \sum_{h \in [m_i]} \ell(x, d_{i,h})
\end{equation}
where $\ell : \R^n \times \R^d \to \R$ is a loss function and $d_{i,h} \in \R^d$ data points stored by agent $i$.

Different algorithm designs have been proposed to solve~\eqref{eq:problem}, with the main categories being \textit{gradient-based} and \textit{dual-based}.
In this article we will offer a unifying perspective on both classes of algorithms, through the lens of non-expansive operator theory. This viewpoint provides interesting insights and can spark further original research.

The article will provide in section~\ref{sec:preliminaries} a short primer on operator theory and its tools, with an intermezzo on average consensus. Operator theory will then be leveraged when discussing distributed algorithms in section~\ref{sec:distributed-algorithms}. Section~\ref{sec:conclusion} will conclude with a discussion of these algorithms in light of the practical challenges that arise when deploying them.

\section{Preliminaries on Operator Theory and Consensus}\label{sec:preliminaries}
In this section we review the necessary background in operator theory which will underlie section~\ref{sec:distributed-algorithms}. We further discuss these background notions in light of \textit{average consensus}, foundational tool in distributed optimization and learning.

\subsection{Background on Operator Theory}\label{subsec:background-operators}
This section presents a limited review of operator theory, geared towards discussing distributed optimization and learning; for a more comprehensive background, we refer to \cite{ryu_primer_2016,bauschke_convex_2017}.

Consider an operator $\T(\cdot) : \R^n \to \R^n : \x \mapsto \T(\x)$.
We say that $\bar{x} \in \R^n$ is a \textit{fixed point} of $\T$ if $\bar{\x} = \T(\bar{\x})$, and we denote by $\fix(\T) = \{ \bar{\x} \in \R^n \ | \ \bar{\x} = \T(\bar{\x}) \}$ the \textit{fixed point set}.
The central goal in operator theory is to describe under what properties of $\T$ its repeated application $\x_{k+1} = \T(\x_k)$, $k \in \N$, converges to a fixed point.
We start by defining the most important of such properties.

\begin{definition}[Non-expansive, contractive operators]\label{def:lipschitz-operators}
Operator $\T$ is \textit{$\zeta$-Lipschitz continuous}, $\zeta \geq 0$, if for all $\x, \y \in \R^n$:
$
	\norm{\T(\x) - \T(\y)} \leq \zeta \norm{\x - \y}.
$
The operator is \textit{non-expansive} if $\zeta = 1$, and \textit{contractive} if $\zeta \in (0, 1)$.
\end{definition}

\begin{definition}[Averaged operators]\label{def:averaged-operators}
Operator $\T$ is \textit{$\alpha$-averaged} if and only if there exist $\alpha \in (0,1)$ and a non-expansive operator $\mathsf{R} : \R^n \to \R^n$ such that $\T = (1-\alpha) \I + \alpha \mathsf{R}$, with $\I$ denoting the identity operator.
Equivalently, $\T$ is $\alpha$-averaged if for all $\x, \y \in \R^n$
$$
	\norm{\T(\x) - \T(\y)}^2 \leq \norm{\x - \y}^2 - \frac{1 - \alpha}{\alpha} \norm{(\I - \T)(\x) - (\I - \T)(\y)}^2.
$$
\end{definition}

Before we can prove convergence to a fixed point, we need to guarantee that $\fix(\T)$ is indeed non-empty, resorting to the following results (cf. Theorems~4.29 and~1.50 in \cite{bauschke_convex_2017}).

\begin{theorem}[Fixed points of non-expansive operators]\label{th:browder}
Let $\D \subset \R^n$ be non-empty, convex, and compact, and assume that $\T : \D \to \D$ is non-expansive; then $\fix(\T) \neq \emptyset$.
\end{theorem}

\begin{theorem}[Fixed points of contractive operators]
Let $\T : \R^n \to \R^n$ be $\zeta$-contractive; then there exists a unique fixed point, $\fix(\T) = \{ \bar{\x} \}$.
\end{theorem}

Clearly, contractiveness is a stronger property than non-expansiveness and averagedness. This is also reflected on their convergence, as exemplified by the following two results.

\begin{theorem}[Convergence of non-expansive operators]\label{th:convergence-nonexpansive}
Let $\D \subset \R^n$ be non-empty, convex, and closed, and assume that $\T : \D \to \D$ is non-expansive with $\fix(\T) \neq \emptyset$ (\textit{e.g.} $\D$ is also bounded and Theorem~\ref{th:browder} applies). Then for any $\alpha \in (0, 1)$ the trajectory $\{ \x_k \}_{k \in \N}$ generated by
\begin{equation}\label{eq:km}
    \x_{k+1} = (1 - \alpha) \x_k + \alpha \T(\x_k), \,\, \x_0 \in \D
\end{equation}
converges to a point in $\fix(\T)$. Additionally, the following bound holds for any $\bar{\x} \in \fix(\T)$:
\begin{equation}\label{eq:fpr-bound}
    \norm{(\I - \T) \x_k} \leq \frac{1}{\sqrt{k+1}} \sqrt{\frac{\alpha}{1 - \alpha}} \norm{\x_0 - \bar{\x}}.
\end{equation}
\end{theorem}

The repeated application of a non-expansive operator, $\x_{k+1} = \T(\x_k)$, does not converge in general, since it simply guarantees that $\norm{\x_{k+1} - \bar{\x}} \leq \norm{\x_k - \bar{\x}} \leq \ldots \leq \norm{\x_0 - \bar{\x}}$. For this reason, Theorem~\ref{th:convergence-nonexpansive} analyzes the repeated application of the \textit{relaxed} operator $(1 - \alpha) \I + \alpha \T$, $\alpha \in (0,1)$, named \km iteration. Convergence to a fixed point of $\T$ is then guaranteed by the fact that $\fix(\T) \equiv \fix((1 - \alpha) \I + \alpha \T)$.

\begin{theorem}[Convergence of contractive operators]\label{th:convergence-contractive}
Let $\T : \R^n \to \R^n$ be $\zeta$-contractive. Then the trajectory $\{ \x_k \}_{k \in \N}$ generated by
\begin{equation}\label{eq:bp}
    \x_{k+1} = \T(\x_k), \,\, \x_0 \in \R^n
\end{equation}
converges to the unique fixed point $\bar{\x}$. Additionally, the following bound holds:
\begin{equation}\label{eq:bp-bound}
    \norm{\x_k - \bar{\x}} \leq \zeta^k \norm{\x_0 - \bar{\x}}.
\end{equation}
\end{theorem}

Theorem~\ref{th:convergence-contractive} thus shows that for a contractive $\T$ we can prove \textit{linear} convergence of $\{ \x_k \}_{k \in \N}$ to the fixed point. On the other hand, by Theorem~\ref{th:convergence-nonexpansive} we see that convergence to a fixed point is \textit{sub-linear} for non-expansive operators -- this follows from~\eqref{eq:fpr-bound} and noting that $\norm{(\I - \T) \bar{\x}} = 0$ if and only if $\bar{\x} \in \fix(\T)$.

\begin{remark}[A dynamical systems perspective]
We can draw a parallel between the operator theory formalism and that of dynamical systems theory. The idea is to interpret the repeated application of an operator $\T$ as the update equation of a discrete-time dynamical system $\x_{k+1} = \T(\x_k)$.
Thus, the fixed points of $\T$ are the \textit{equilibria} of the system, whose (asymptotic) stability depends on the properties of the operator.
For example, if $\T$ is contractive, then $\bar{\x}$ is \textit{(globally) asymptotically stable}, with $V(\x) = \norm{\x - \bar{\x}}^2$ serving as a Lyapunov function. On the other hand, the fixed points of a non-expansive $\T$ are only \textit{stable}.
\end{remark}

\subsection{Intermezzo: Average Consensus}\label{subsec:average-consensus}
In this section we review average consensus through the lens of operator theory. Consider a multi-agent system modeled by the graph $\mathcal{G} = (\mathcal{V}, \Ecal)$, $\mathcal{V} = [N]$; the goal is to enable the agents to compute the average of their local observations $\{ u_i \}_{i \in [N]}$, but only resorting to peer-to-peer communications.
The most widely studied average consensus protocol is characterized by
\begin{equation}\label{eq:average-consensus}
    x_{i,k+1} = w_{ii} x_{i,k} + \sum_{j \in \Ncal_i} w_{ij} x_{j,k}, \quad x_{i,0} = u_i, \quad i \in [N], \quad k \in \N,
\end{equation}
where $\Ncal_i$ denotes the neighbor set of agent $i$, and $w_{ij} \geq 0$. Clearly, agent $i$ can update its state $x_i$ only after receiving its neighbors' states.
We can interpret~\eqref{eq:average-consensus} as the application of a \textit{linear} operator
$
    \x_{k+1} = \W \x_k
$
where $\x \in \R^N$ collects the local states $x_i$, and $\W \in \R^{N \times N}$ collects the consensus weights $w_{ij}$. Assuming $\W$ to be row stochastic, we can see that the fixed set of $\W$ is the consensus set $\fix(\W) = \{ \x \in \R^N \ | \ x_1 = \ldots = x_N \}$.
The consensus protocol can also be interpreted as the gradient descent (cf. Definition~\ref{def:gradient-descent}) applied to $\min \frac{1}{2} \x^\top (\Ib - \W) \x$.

\begin{figure}[!ht]
\begin{subfigure}{0.3\textwidth}
\centering
\includegraphics*[scale=0.8]{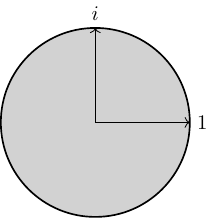}
\caption{Non-expansive.}
\end{subfigure}
\begin{subfigure}{0.3\textwidth}
\centering
\includegraphics*[scale=0.8]{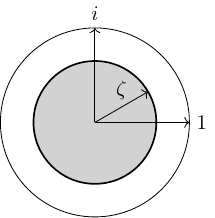}
\caption{$\zeta$-contractive.}
\end{subfigure}
\begin{subfigure}{0.3\textwidth}
\centering
\includegraphics*[scale=0.8]{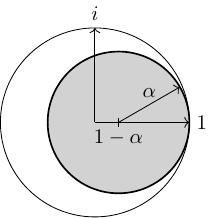}
\caption{$\alpha$-averaged.}
\label{fig:linear-operators-avg}
\end{subfigure}
\caption{Depiction of where the eigenvalues of linear operators lie, depending on the operators' properties.}
\label{fig:linear-operators}
\end{figure}

Now, the operator theoretical properties of Definitions~\ref{def:lipschitz-operators} and~\ref{def:averaged-operators} have an interesting interpretation for linear operators in terms of where their eigenvalues lie (\cite{iutzeler_generic_2019}), as depicted in Figure~\ref{fig:linear-operators}.
Inspecting the figures, we see that non-expansiveness is not sufficient to ensure convergence, as the eigenvalues can lie anywhere on the unit circle. Moreover, for an $\alpha$-averaged operator the eigenvalues lie in a disk with center $1 - \alpha$ and radius $\alpha$, with $1$ being a semi-simple eigenvalue.

Further assuming that $\W$ is symmetric (and thus also column stochastic), its eigenvalues satisfy $-1 < \lambda_1 \leq \ldots \leq \lambda_N = 1$, and hence $\W$ is $(1 - \lambda_1) / 2$-averaged by Figure~\ref{fig:linear-operators-avg}. Theorem~\ref{th:convergence-nonexpansive} then ensures convergence to one of the fixed points, and in particular to the average consensus $\lim_{k \to \infty} x_{i,k} = \frac{1}{N} \sum_{i \in [N]} u_i$ by column stochasticity and the initialization $x_{i,0} = u_i$ (see \textit{e.g.} \cite{olfati_consensus_2007}).
We remark that averaged linear operators converge linearly, whereas averaged non-linear operators converge only sub-linearly (cf. Theorem~\ref{th:convergence-nonexpansive}), as opposed to contractive (non-linear) operators which converge linearly (cf. Theorem~\ref{th:convergence-contractive}).

\subsection{Operator Theory for Convex Optimization}
In this section we review the application of operator theory for convex optimization, presenting some of the most common algorithms. In the following, we denote by $\ccp(\R^n)$ the set of closed, convex and proper functions
\footnote{A function $f$ is proper if its domain is non-empty and it takes values in $(-\infty, +\infty]$; closed if its epigraph $\{ (x, t) \in \R^n \times \R \ | \ f(x) \leq t \}$ is a closed set.}
We start by defining three fundamental operators.

\begin{figure}[!ht]
\centering
    \begin{subfigure}{0.4\textwidth}
    \includegraphics{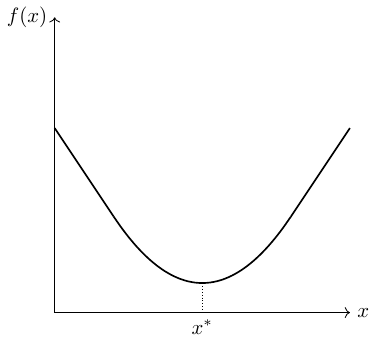}
    \end{subfigure}

    \begin{subfigure}{0.4\textwidth}
    \includegraphics{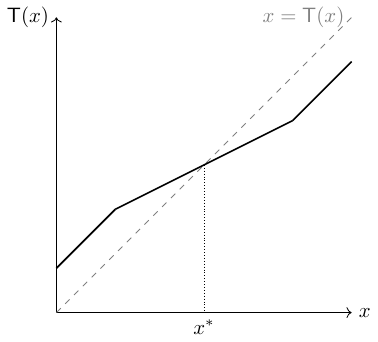}
    \end{subfigure}%
\caption{Huber loss (top) and gradient descent operator applied to it (bottom), highlighting that the fixed point coincides with the minimizer.}
\label{fig:minima-fixed-points}
\end{figure}

\begin{definition}[Gradient descent operator]\label{def:gradient-descent}
Let $f \in \ccp(\R^n)$ be $\lmax$-smooth, we define the \textit{gradient descent operator} as $\I - \rho \nabla f(\cdot)$, where $\rho > 0$ is the step-size.
The fixed points of the gradient descent operator coincide with the minimizers of $f$\footnote{This is a simple consequence of Fermat's rule requiring that the minimizers of $f$ satisfy $\nabla f(\x) = 0$.}, as exemplified by Figure~\ref{fig:minima-fixed-points}.
The operator is averaged if $\rho < 2 / \lmax$, and contractive if additionally $f$ is $\lmin$-strongly convex (\cite{ryu_primer_2016}).
\end{definition}

\begin{definition}[Proximal and reflective operators]
Let $f \in \ccp(\R^n)$, we define the \textit{proximal operator} of $f$ as
$$
    \prox_{\rho f}(\y) = \argmin_{\x \in \R^n} \left\{ f(\x) + \frac{1}{2 \rho} \norm{\x - \y}^2 \right\}
$$
where $\rho > 0$ is the \textit{penalty parameter}. Moreover, we define the \textit{reflective operator} as
$$
    \refl_{\rho f}(\cdot) = 2 \prox_{\rho f}(\cdot) - \I.
$$
The fixed points of the proximal and reflective operators coincide with the minimizers of $f$.
The $\prox_{\rho f}$ is $1/2$-averaged, and $\refl_{\rho f}$ is non-expansive (\cite{bauschke_convex_2017}); moreover, if $f$ is $\lmin$-strongly convex then they are both contractive (\cite{giselsson_linear_2017}).
\end{definition}

We can see that convexity results in averagedness of the gradient and proximal operators, while the more powerful assumption of strong convexity ensures contractiveness. In turn this reflects on the convergence properties of the operators, as observed in section~\ref{subsec:background-operators}.
Now, these operators are important in their own right, especially the gradient descent operator; however, we are interested in their use as building blocks of the more complex operators that we will need in section~\ref{sec:distributed-algorithms}.
In particular, consider the \textit{composite} problem $\min_{\x \in \R^n} f(\x) + g(\x)$, then the following operators can be applied to solve it.

\begin{definition}[Proximal gradient operator]\label{def:proximal-gradient}
Let $f, g \in \ccp(\R^n)$, and assume that $f$ is $\lmax$-smooth. We define the \textit{proximal gradient} operator as
$$
    \prox_{\rho g} \comp \left( \I - \rho \nabla f(\cdot) \right)
$$
with $\rho \in (0, 2 / \lmax)$.
The fixed points coincide with the minimizers of $f + g$; the operator is averaged, and it is contractive if additionally $f$ is $\lmin$-strongly convex (\cite{ryu_primer_2016}).
\end{definition}

\begin{definition}[Peaceman-Rachford operator]\label{def:prs}
Let $f, g \in \ccp(\R^n)$, we define the \textit{Peaceman-Rachford} operator as
$$
    \refl_{\rho g} \comp \refl_{\rho f}
$$
with $\rho > 0 $.
The fixed points coincide with the minimizers of $f + g$; the operator is non-expansive, and it is contractive if $f$ is $\lmax$-smooth and $\lmin$-strongly convex (\cite{giselsson_linear_2017}).
\end{definition}

The proximal gradient and Peaceman-Rachford operators are called \textit{splitting}, as they apply a single operation (\textit{e.g.} gradient or reflective) to $f$ and $g$ separately.
We remark that, in the absence of strong convexity, the Peaceman-Rachford is only non-expansive, and thus its \textit{relaxed} version $(1 - \alpha) \I + \alpha \refl_{\rho g} \comp \refl_{\rho f}$ is usually applied in practice. Relaxation ensures that the overall operator is averaged, and by the fact that $\fix(\T) \equiv \fix((1 - \alpha) \I + \alpha \T)$, Theorem~\ref{th:convergence-nonexpansive} thus ensures convergence to a fixed point of the Peaceman-Rachford operator.

\section{Distributed Optimization and Learning}\label{sec:distributed-algorithms}
In this section we discuss different distributed optimization algorithms in an operator theoretical perspective. We focus on the two most studied classes of distributed algorithms: gradient-based and ADMM. The final section will discuss and compare these alternative approaches. For a comprehensive overview of distributed optimization we reference \cite{sayed_adaptation_2014,nedic_distributed_2018,notarstefano_distributed_2019}.

\subsection{Distributed Gradient Methods}
We start by noting that the consensus optimization problem~\eqref{eq:problem} can be written as $\min_{\x \in \R^{nN}} f(\x) + g(\x)$, where $\x$ collects the local states $x_i \in \R^n$, $i \in [N]$, $f(\x) = \sum_{i \in [N]} f_i(x_i)$ and
$$
    g(\x) = \begin{cases}
        0 & \text{if} \ x_1 = \ldots = x_N, \\
        +\infty & \text{otherwise}.
    \end{cases}
$$
Thus, in principle we could apply the proximal gradient of Definition~\ref{def:proximal-gradient} to solve it, which would yield the updates
\begin{equation}\label{eq:projected-gradient}
    x_{i,k+1} = \frac{1}{N} \sum_{j \in [N]} \left( x_{j,k} - \rho \nabla f_j(x_{j,k}) \right), \quad i \in [N], \quad k \in \N,
\end{equation}
since \textit{the proximal of $g$ corresponds to a projection onto the consensus set}, that is, to the average of the components of $\x_k - \rho \nabla f(\x_k)$.
However, computing such projection would require communications between each agent and all others, which is not possible in a distributed set-up where each agent is connected to only a subset of the others. It is therefore necessary to design a distributed projection onto the consensus set.

\subsubsection{Distributed gradient descent}
The first approach one can apply is to replace the projection with average consensus, reviewed in section~\ref{subsec:average-consensus}. Indeed, average consensus can be performed in a distributed fashion, relying only on peer-to-peer communications.
This results in the algorithm characterized by the following updates:
\begin{equation}\label{eq:dgd}
    x_{i,k+1} = \sum_{j \in \Ncal_i \cup \{ i \}} w_{ij} \left( x_{j,k} - \rho \nabla f_j(x_{j,k}) \right), \quad i \in [N], \quad k \in \N,
\end{equation}
which can be represented as the operator update
$$
    \x_{k+1} = (\W \otimes I_n) \comp \left( \I - \rho \nabla f(\cdot) \right) (\x_k),
$$
where $\otimes$ denotes Kronecker product.
This algorithm is called Adapt-Then-Combine (\cite{chen_distributed_2013}), and another widely studied variant is that of DGD, defined by $\x_{k+1} = (\W \otimes I_n) \x_k - \rho \nabla f(\x_k)$ (\cite{yuan_convergence_2016}).

However, average consensus is a poor substitute for a full projection onto the consensus set. Indeed, this class of algorithms has a fundamental issue: \textit{their fixed points do not coincide with solutions of~\eqref{eq:problem}}.
In particular, it can be proved that a fixed point $\bar{\x}$ of~\eqref{eq:dgd} lies in a ball centered in $\x^*$ with radius proportional to the step-size $\rho$, see \textit{e.g.} \cite[Theorem~3]{chen_distributed_2013}.
Therefore,~\eqref{eq:dgd} is guaranteed to converge by the results reviewed in section~\ref{sec:preliminaries} -- but to the wrong point.
To ensure exact convergence of~\eqref{eq:dgd} it is necessary to employ a diminishing step-size $\rho_k$, chosen so that $\sum_{k \in \N} \rho_k = +\infty$, $\sum_{k \in \N} \rho_k^2 < +\infty$ (\cite{nedic_distributed_2018}). However, convergence is attained with only a sub-linear rate, even for strongly convex functions.

\subsubsection{Gradient tracking}\label{subsec:gradient-tracking}
The previous section showed that approximating a projection onto the consensus set with one step of average consensus is not sufficient to guarantee convergence to a solution $\x^*$ of~\eqref{eq:problem}. Indeed, the agents perform consensus on the components of the vector $\x_k - \rho \nabla f(\x_k)$, which changes at each iteration, and average consensus is unable to track the average of a time-varying signal.
The idea then is to employ a distributed protocol that \textit{can} track a time-varying average: \textit{dynamic average consensus} (\cite{kia_tutorial_2019}). Let $\{ y_{i,k} \}_{k \in \N}$, $i \in [N]$, be time-varying signals observed by the agents, then dynamic average consensus is characterized by
\begin{equation}\label{eq:dynamic-consensus}
    x_{i,k+1} = w_{ii} x_{i,k} + \sum_{j \in \Ncal_i} w_{ij} x_{j,k} + y_{i,k} - y_{i,k-1}, \quad i \in [N], \quad k \in \N,
\end{equation}
or, in vector form $\x_{k+1} = (\W \otimes I_n) \x_k + \y_k - \y_{k-1}$. Notice that~\eqref{eq:dynamic-consensus} can be interpreted as an affine operator with a time-varying offset.

The idea now is to apply dynamic consensus~\eqref{eq:dynamic-consensus} to approximate the consensus projection of~\eqref{eq:projected-gradient}, which yields the \textit{gradient tracking} algorithm
\begin{equation}\label{eq:gt}
\begin{split}
    y_{i,k} &= x_{i,k} - \rho \nabla f_i(x_{i,k}) \\
    x_{i,k+1} &= w_{ii} x_{i,k} + \sum_{j \in \Ncal_i} w_{ij} x_{j,k} + y_{i,k} - y_{i,k-1}
\end{split} \qquad  i \in [N], \quad k \in \N.
\end{equation}
The use of dynamic consensus guarantees that the algorithm converges to an optimal solution of~\eqref{eq:problem}; or, in other words, the fixed points of~\eqref{eq:gt} coincide with the solutions of~\eqref{eq:problem}.
Different variants of gradient tracking algorithms have been proposed, and we reference \cite{bof_multiagent_2019,jakovetic_unification_2019,xin_general_2020,xu_distributed_2021} for a comprehensive review. These algorithms differ in \textit{(i)} where dynamic average consensus is applied, \textit{e.g.} to the states, to the gradients, or to the gradient steps (as in~\eqref{eq:gt}), and \textit{(ii)} if additional consensus rounds are employed, to improve performance.

Interestingly, gradient tracking algorithms admit a second interpretation as \textit{primal-dual operators} (\cite{alghunaim_linear_2020}). To see this, the first step is to define the following equivalent form of~\eqref{eq:problem}
\begin{equation}\label{eq:problem-gt}
    \x^* \in \argmin_{x_i \in \R^n, \ i \in [N]} \sum_{i \in [N]} f_i(x_i) \quad \text{s.t.} \quad (\Ib - \W )^{1/2} \otimes I_n \ \x = \0
\end{equation}
where $\W$ is a consensus matrix, and further to rewrite it as the \textit{saddle point} problem
\begin{equation}\label{eq:problem-saddle}
    \min_{\x} \max_{\w} \Lcal(\x, \w) := f(\x) + \langle \w, (\Ib - \W )^{1/2} \otimes I_n \x \rangle
\end{equation}
where $\Lcal$ is the Lagrangian of~\eqref{eq:problem-gt}, and $\w \in \R^{n N}$ is the vector of dual variables.
We apply now the primal-dual operator
\begin{subequations}\label{eq:primal-dual}
\begin{align}
    \x_{k+1} &= \x_k - \rho \nabla_{\x} \Lcal(\x_k, \w_k) = \x_k - \rho \left( \nabla f(\x_k) + (\Ib - \W )^{1/2} \otimes I_n \w_k \right) \label{eq:primal-dual-x} \\
    \w_{k+1} &= \w_k + \frac{1}{\rho} \nabla_{\w} \Lcal(\x_{k+1}, \w_k) = \w_k + \frac{1}{\rho} (\Ib - \W )^{1/2} \otimes I_n \x_{k+1} \label{eq:primal-dual-w}
\end{align}
\end{subequations}
where \eqref{eq:primal-dual-x} is the gradient descent operator in $\x$ applied to $\Lcal$, with step-size $\rho$, and \eqref{eq:primal-dual-w} is the gradient ascent operator in $\w$ applied to $\Lcal$, with step-size $1/\rho$. Notice that this is an \textit{incremental} scheme, as \eqref{eq:primal-dual-w} employs the output $\x_{k+1}$ of~\eqref{eq:primal-dual-x}.
Defining $\y_k = \x_k - \nabla f(\x_k)$ and rearranging~\eqref{eq:primal-dual} then yields exactly~\eqref{eq:gt}.

\subsection{Distributed ADMM}
In this section we discuss the main alternative to gradient-based methods: the distributed ADMM (Alternating Direction Method of Multipliers).

\subsubsection{Problem reformulation}
The first step to design the distributed ADMM is to reformulate problem~\eqref{eq:problem} into (\cite{bastianello_asynchronous_2021})
\begin{equation}\label{eq:problem-admm}
\begin{split}
    &\min \sum_{i \in [N]} f_i(x_i) \\
    &\text{s.t.} \ x_i = y_{ij}, \,\, x_j = y_{ji}, \,\, y_{ij} = y_{ji} \quad \forall (i,j) \in \Ecal
\end{split}
\end{equation}
by replacing the consensus constraints with edge-based constraints that employ the auxiliary variables $y_{ij}$. The idea is to enforce agreement between any two neighboring agents, which -- by connectedness of the graph -- is equivalent to ensuring network-wide consensus.
Defining a new cost function $g : \R^{2 n |\Ecal|} \to \R \cup \{ +\infty \}$ such that $g(\y) = 0$ if $y_{ij} = y_{ji}$ $\forall (i,j) \in \Ecal$, $+\infty$ otherwise, and the matrix $\A = \blkdiag\{ \1_{|\Ncal_i|} \}$, we can rewrite~\eqref{eq:problem-admm} as
\begin{equation}\label{eq:problem-admm-2}
    \min \sum_{i \in [N]} f_i(x_i) + g(\y) \quad \text{s.t.} \ \A \x = \y.
\end{equation}
Notice that the particular structure of $\A$ derives from the fact that each agent $i$ participates in the $|\Ncal_i|$ constraints $x_i = y_{ij}$, $j \in \Ncal_i$ -- this is thus the \textit{arc incidence} matrix.

Problem~\eqref{eq:problem-admm-2} is a \textit{composite problem with a linear constraint}, which by strong duality can be reduced to an unconstrained composite problem by switching to its dual \cite[sec.~5.2.3]{boyd_convex_2004}
\begin{equation}\label{eq:problem-admm-3}
    \min_{\w \in \R^{2 n |\Ecal|}} d_f(\w) + d_g(\w)
\end{equation}
where $d_f(\w) = f^*(\A^\top \w)$ and $d_g(\w) = g^*(- \w)$, with the convex conjugate being defined as $f^*(\z) = \sup_{\x} \{ \langle \z, \x \rangle - f(\x) \}$ \cite[sec.~5.1.6]{boyd_convex_2004}.

\subsubsection{Algorithm design}
The dual problem~\eqref{eq:problem-admm-3} is composite, which means that \textit{we can apply the (relaxed) Peaceman-Rachford operator of Definition~\ref{def:prs} to it}, yielding the update $\z_{k+1} = (1 - \alpha) \z_k + \alpha \refl_{\rho d_g} \circ \refl_{\rho d_f}(\z_k)$, $\alpha \in (0, 1)$, $\rho > 0$.
Owing to the \textit{splitting} structure of the operator and the fact that $\refl_{\rho f} = 2 \prox_{\rho f} - \I$, we can equivalently characterize the algorithm by the updates
\begin{subequations}\label{eq:prs}
\begin{align}
    \w_{k+1} &= \prox_{\rho d_f}(\z_k) \label{eq:prs-1} \\
    \uv_{k+1} &= \prox_{\rho d_g}(2 \w_{k+1} - \z_k) \label{eq:prs-2} \\
    \z_{k+1} &= \z_k + 2 \alpha (\uv_{k+1} - \w_{k+1}). \label{eq:prs-3}
\end{align}
\end{subequations}
To summarize, designing the distributed ADMM requires reformulating the consensus constraints of~\eqref{eq:problem} with the edge-based constraints of~\eqref{eq:problem-admm}, yielding a composite problem with linear constraints. The dual of problem~\eqref{eq:problem-admm} is then an unconstrained composite problem, which allows us to apply the Peaceman-Rachford operator.

The updates of~\eqref{eq:prs}, however, reference only the dual (and auxiliary) variables. The last step leading to the distributed ADMM then is to provide a way of computing the proximals in~\eqref{eq:prs-1} and~\eqref{eq:prs-2} in a distributed fashion, which will produce the primal updates. This is possible by using the fact that (\cite{davis_convergence_2016})
\begin{align*}
    \x_{k+1} &= \argmin_{\x} \left\{ f(\x) - \langle \z_k, \A \x \rangle + \frac{\rho}{2} \norm{\A \x}^2 \right\} \\
    \w_{k+1} &= \z_k - \rho \A \x_{k+1} = \prox_{\rho d_f}(\z_k),
\end{align*}
which shows that computing $\prox_{\rho d_f}$ can be done by solving a minimization in $f$. Similar formulae apply to the computation of $\prox_{\rho d_g}$.
Now, leveraging the specific structure of the distributed problem~\eqref{eq:problem-admm}, namely that $f(\x) = \sum_{i \in [N]} f_i(x_i)$, $\A = \blkdiag\{ \1_{|\Ncal_i|} \}$, and the definition of $g(\y)$, we can simplify~\eqref{eq:prs} to (\cite{bastianello_asynchronous_2021}):
\begin{subequations}\label{eq:admm}
\begin{align}
    x_{i,k+1} &= \prox_{f_i / (\rho |\Ncal_i|)}\left( \sum\nolimits_{j \in \Ncal_i} z_{ij,k} / (\rho |\Ncal_i|) \right) \label{eq:admm-1} \\
    z_{ij,k+1} &= (1 - \alpha) z_{ij,k} - \alpha \left( z_{ji,k} - 2 \rho x_{j,k+1} \right). \label{eq:admm-2}
\end{align}
\end{subequations}
The updates~\eqref{eq:admm} fully characterize the distributed ADMM, since they can be performed in a distributed fashion. In particular, each agent $i \in [N]$ stores and updates the variables $x_i$ and $\{ z_{ij} \}_{j \in \Ncal_i}$, and it only needs to receive the vectors $\{ z_{ji,k} - 2 \rho x_{j,k+1} \}_{j \in \Ncal_i}$ through peer-to-peer communications.

\subsubsection{Alternative interpretation}
The discussion above shows how the distributed ADMM can be derived by applying the Peaceman-Rachford operator to the dual of the (suitably reformulated) consensus optimization problem.
However, there is a second interpretation of this algorithm: \textit{as a primal-dual method}, similarly to the gradient tracking in section~\ref{subsec:gradient-tracking}.
To see this, let us define the \textit{augmented Lagrangian} of~\eqref{eq:problem-admm-2} as
$$
    \Lcal(\x, \y, \w) = f(\x) + g(\y) + \langle \w, \A \x - \y \rangle + \frac{\rho}{2} \norm{\A \x - \y}^2;
$$
notice that this Lagrangian differs from that of~\eqref{eq:problem-saddle} due to the regularization (augmentation) term.
Applying alternating minimization in the primal variables $\x$, $\y$, and a linear update in the dual $\w$ yields the algorithm
\begin{equation}\label{eq:lagrangian-admm}
\begin{split}
    \x_{k+1} &= \argmin_{\x} \Lcal(\x, \y_k, \w_k) \\
    \y_{k+1} &= \argmin_{\y} \Lcal(\x_{k+1}, \y, \w_k) \\
    \w_{k+1} &= \w_k + \rho \left( \A \x_{k+1} - \y_{k+1} \right).
\end{split}
\end{equation}
Exploiting the particular distributed structure of~\eqref{eq:problem-admm-2} it is then possible to show the equivalence of the Lagrangian ADMM~\eqref{eq:lagrangian-admm} and~\eqref{eq:admm} (with $\alpha = 1/2$, but the equivalence can be generalized to any $\alpha$) (\cite{bastianello_asynchronous_2021}).

\smallskip

We conclude this section by summarizing in Table~\ref{tab:summary-interpretations} the operator theoretical interpretations of distributed algorithms reviewed above.
\begin{table}[!ht]
\TBL{\caption{Summary of operator theoretical interpretations of distributed algorithms}\label{tab:summary-interpretations}}
{\begin{tabular*}{\textwidth}{ll}
\toprule
\TCH{Algorithm} &
\TCH{Interpretation} \\
\colrule
Consensus & Linear, averaged operator \\
& Gradient descent on $\frac{1}{2} \x^\top (\Ib - \W) \x$\\
Distributed gradient descent & Projected gradient with $\proj \approx$ average consensus \\
Gradient tracking & Projected gradient with $\proj \approx$ \textit{dynamic} average consensus \\
 & \& Primal-dual \\
ADMM & Peaceman-Rachford applied to dual problem \\
 & \& Primal-dual (with augmented Lagrangian) \\
\botrule
\end{tabular*}}{}
\end{table}

\section{Conclusion: The Challenges of Distributed Optimization and Learning}\label{sec:conclusion}
We conclude by comparing the algorithms presented in section~\ref{sec:distributed-algorithms}, and discussing their features in light of the practical challenges that arise when deploying them, see \textit{e.g.} \cite{li_federated_2020}. The final section will discuss some recent trends in both distributed and federated optimization and learning.

\paragraph{Discussion}
The first comparison we draw is on the number of local variables (intended as vectors in $\R^n$) that each agent needs to store while executing a specific distributed algorithm. As summarized in Table~\ref{tab:comparison}, distributed gradient descent has the smallest memory footprint, as the agents need to store a single vector each. However, this is demonstrably not enough to achieve convergence, since all distributed gradient methods with a single local state fail to converge \cite[Corollay~6]{sundararajan_canonical_2019}; this result complements the aforementioned \cite[Theorem~3]{chen_distributed_2013}.
\begin{table}[!ht]
\TBL{\caption{Comparison of the distributed algorithms in section~\ref{sec:distributed-algorithms}}\label{tab:comparison}}
{\begin{tabular*}{\textwidth}{lllllll}
\toprule
\TCH{Algorithm} &
\TCH{\# local var.s} &
\TCH{Exact conv.} &
\TCH{Asynchrony} &
\TCH{Lossy comm.} &
\TCH{Noisy comm.}\footnotemark{a} &
\TCH{Local computations}\\
\colrule
Distributed gradient descent & $1$ & \xmark & \cmark & \cmark & \cmark & $\I - \rho \nabla f_i(\cdot)$ \\
Gradient tracking & $2$ & \cmark & \xmark & \xmark & \xmark & $\I - \rho \nabla f_i(\cdot)$ \\
Robust gradient tracking & $2|\Ncal_i| + 4$ & \cmark & \cmark & \cmark & \xmark & $\I - \rho \nabla f_i(\cdot)$ \\
ADMM & $|\Ncal_i| + 1$ & \cmark & \cmark & \cmark & \cmark & $\prox_{f_i}$ \\
\botrule
\end{tabular*}}{%
\begin{tablenotes}
\footnotetext[a]{\textit{e.g.} due to quantization}
\end{tablenotes}
}%
\end{table}
Algorithms that store multiple local states are then necessary, as the gradient tracking algorithm~\eqref{eq:gt} which requires two states (in the formulation of~\eqref{eq:primal-dual}) in order to guarantee exact convergence.
ADMM provides the same convergence guarantee, but with a larger storage requirement. Indeed, the number of local variables scales with the number of peer-to-peer connections, since the constraints of~\eqref{eq:problem-admm} introduce two auxiliary variables for each edge.
However, as discussed in the following, this feature of ADMM is actually the key to its robustness in the face of different practical challenges.

The agents collaborating towards the solution of~\eqref{eq:problem} are in general equipped with \textit{heterogeneous resources}, for example in terms of computational power, memory, battery, \textit{etc.} As a consequence, different agents will complete local computations (\textit{e.g.} gradient evaluations) with different speed and accuracy, resulting in \textit{asynchrony}.
ADMM has been proved to be robust to asynchrony (\cite{bastianello_asynchronous_2021}), while gradient tracking is not, due to its reliance on the synchronous dynamic consensus in~\eqref{eq:gt}. Indeed, to ensure robustness of gradient tracking, the distributed implementation of the consensus projection needs to be replaced with a robust version (\cite{bof_multiagent_2019,tian_achieving_2020}). Importantly, this robust consensus also requires a number of additional variables that scales with the number of edges.
From an operator theoretical perspective, an asynchronous algorithm can be modeled as a \textit{random coordinate update operator} (\cite{peng_coordinate_2016}), where for example coordinate $i$ is updated only when agent $i$ activates.

Another repercussion of resource limitations is that the agents may need to resort to inexact local computations.
For example, in ADMM the local proximal updates~\eqref{eq:admm-1} entail the solution of a minimization problem, which in practice may need to be approximated. Similarly, for gradient-based methods the evaluation of a full gradient may be computationally expensive when the empirical risk minimization costs~\eqref{eq:erm} are defined on a large dataset. The agents may then resort to stochastic gradients instead.
The use of inexact local computations (\textit{e.g.} approximate proximals or stochastic gradients) results in \textit{inexact operators} of the type $\x_{k+1} = \T(\x_k) + \e_k$, for some computation error $\e_k$.

Alongside the local resource limitations, the agents are also subject to network-wide \textit{communications constraints}. Indeed, distributed algorithms rely on peer-to-peer communications which in practice may have limited bandwidth, especially when wireless connections are employed.
Thus the number and size of communications that can be exchanged at each iteration $k$ may be constrained, for example through the use of quantization, which similarly to inexact local computations yields additive operator errors.
Additionally, communications may get lost due to interference, requiring thus robust algorithm designs (\cite{bof_multiagent_2019,tian_achieving_2020,bastianello_asynchronous_2021}).

\paragraph{Outlook}
We conclude with a summary of current trends at the intersection of operator theory with distributed and federated optimization and learning.
First of all, we remark that the operator theoretical tools and concepts employed in section~\ref{sec:distributed-algorithms} can similarly be applied in federated scenarios (cf. Figure~\ref{fig:architecture-federated}), both to analyze existing algorithm and to design novel ones (\cite{tran_dinh_feddr_2021,grudzien_can_2023}).

Moreover, operator theory can be used to design algorithms that are robust to asynchrony and packet loss, as well as to reduce their communication footprint. For example, it has been applied for compression, to reduce the size of communications (\cite{takezawa_communication_2023}), as well as to reduce their frequency (\cite{grudzien_can_2023}).

Finally, other objectives that we envision to be addressed with operator theory are privacy, robustness to attacks, and personalization (\cite{li_federated_2020}).

\begin{ack}
The work of N. Bastianello was partially supported by the European Union’s Horizon 2020 research and innovation programme under grant agreement No. 101070162.
\end{ack}

\seealso{article title article title}

\bibliographystyle{Harvard}
\bibliography{references}

\end{document}